\documentclass[12pt]{amsart}

\usepackage[english]{babel}
\usepackage{graphicx}
\usepackage{textcomp}
\usepackage{amsmath}

\usepackage{amssymb}
\usepackage{latexsym}
\title[On effects of stochastic
regularization] {On effects of stochastic regularization for the
pressureless gas dynamics}

\author[Korshunova,\, Rozanova]{Anastasia Korshunova $^{1}$,\, Olga Rozanova $^{1}$}

\address[$^{1}$]{Mathematics and Mechanics Faculty, Moscow State University, Moscow
119992, Russia}

\thanks {Supported by  DFG 436 RUS 113/823/0-1 (O.R.) and the special program of the
Ministry of Education of the Russian Federation "The development of
scientific potential of the Higher School", project 2.1.1/1399 (A.K.
and O.R.).}

\email[$^{1}$]{$korshunova{\underline{\phantom{a}}}aa@mail.ru$}
\email[$^{2}$]{\it rozanova@mech.math.msu.su}



\begin {document}

\maketitle

\newtheorem{thm}{Theorem}[section]

\newtheorem{prop}[thm]{Proposition}
\theoremstyle{definition}
\newtheorem{defn}[thm]{Definition}
\theoremstyle{remark}
\newtheorem{rem}[thm]{Remark}
\numberwithin{equation}{section}


\begin{abstract}

We extend our result of \cite{AKR} and show that one can associate
with the stochastically perturbed non-viscid Burgers equation a
system of viscous balance laws. The Cauchy data for the Burgers
equation generates the data for  this system. Till the moment of the
shock formation in the solution to the Burgers equation the above
system of viscous balance laws can be reduced to the pressureless
gas dynamics system (in a limit as the parameters of perturbation
tend to zero). If the solution to the Burgers equation contains
shocks, the limit system is equivalent to the system with a specific
pressure, in some sense analogous to the pressure of barotropic
monoatomic gas.
\end{abstract}


\section{Stochastic perturbation of the Burgers\\ equation}

Let us consider the Cauchy problem for the non-viscous Burgers
equation:
\begin{equation}
\label{equ_Burg} \partial_t u+(u,\nabla)u=0,\,t>0,\qquad
u(x,0)=u_0(x)\in C^1(\mathbb{R}^n)\cap C_b(\mathbb{R}^n).
\end{equation}
Here $u(x,t)=(u_1,...,u_n)(x,t)\,$ is a vector-function
$\mathbb{R}^{n+1}\rightarrow\mathbb{R}^n.$

We associate with (\ref{equ_Burg}) the following  system of
stochastic differential equations:
\begin{equation}
\label{SDU}dX_k(t)=U_k(t)dt+\sigma_1 d(W_k^1)_t,\quad
dU_k(t)=\sigma_2 d(W_k^2)_t,\quad k=1,...,n,\,
\end{equation}
$$X(0)=x,\quad U(0)=u,\quad t>0,$$
where $X(t)$ and $U(t)$ are considered as  random variables with
given initial distributions, $(X(t),U(t))$ runs in the phase space
$\mathbb{R}^n\times\mathbb{R}^n,$  $\sigma_1$ and $\sigma_2$ are
nonnegative constants such that $|\sigma|\ne 0$
($\sigma=(\sigma_1,\sigma_2)$) and $(W^i)_t=(W^i_1,\dots,W^i_n)_t$,
$i=1,2$, are independent $n$ - dimensional Brownian motions.

Let  $P(t,x,u)$ be  a probability density of joint distribution of
random values $(X,U), $ subject to initial data
\begin{equation}
\label{P0}
P_0(x,u)=\delta(u-u_0(x))\rho_0(x)=\prod\limits_{k=1}^n\delta(u_k-(u_0(x))_k)\rho_0(x),
\end{equation}
where $\rho_0(x)$ is a bounded nonnegative
function from $L^1_{loc}({\mathbb R}^n)$.

Let us introduce functions
\begin{equation}
\label{u_sdu} \rho_\sigma(t,x) =
\int\limits_{\mathbb{R}^n}P(t,x,u)\,du, \quad
u_{\sigma}(t,x)=\dfrac{\int\limits_{\mathbb{R}^n}uP(t,x,u)du}{\int\limits_{\mathbb{R}^n}P(t,x,u)du}.
\end{equation}
Notice that $u_{\sigma}(0,x)=u_0(x)$. Certain properties of
$u_{\sigma}(t,x)$ have been established in \cite{ARMMM} and
\cite{ARAMS}.

The density $P=P(t,x,u)$ obeys the Fokker-Planck equation
\begin{equation}
\label{Fok-Plank} \dfrac{\partial P}{\partial
t}=\left[-\sum\limits_{k=1}^n u_k\dfrac{\partial}{\partial
x_k}+\sum\limits_{k=1}^n\dfrac12\sigma_1^2\dfrac{\partial^2}{\partial
x_k^2}+\sum\limits_{k=1}^n\dfrac12\sigma_2^2\dfrac{\partial^2}{\partial
u_k^2}\right] P,
\end{equation}
subject to  initial data (\ref{P0}).

We apply the Fourier transform with respect to  variables $x$ and
$u$  in (\ref{Fok-Plank}), (\ref{P0})   and obtain the following
Cauchy problem for the Fourier transform
$\tilde{P}=\tilde{P}(t,\lambda,\xi)$ of $P(t,x,u)$:
\begin{equation}
\label{preobr_Fok-Plank} \dfrac{\partial \tilde{P}}{\partial
t}=-\dfrac12(\sigma^2_1|\lambda|^2+\sigma^2_2|\xi|^2)\tilde{P}+(\lambda,\dfrac{\partial\tilde{P}}{\partial
\xi}),
\end{equation}
\begin{equation}
\label{preobr_P0}\tilde{P}(0,\lambda,\xi)=\int\limits_{\mathbb{R}^n}e^{-i(\lambda,s)}e^{-i(\xi,u_0(s))}\rho_0(s)ds,
\qquad \lambda,\xi \in {\mathbb R}^n.
\end{equation}
Equation (\ref{preobr_Fok-Plank}) can easily be  integrated and we
obtain the solution:
\begin{equation}
\label{preobr_P}\tilde{P}(t,\lambda,\xi)=\tilde{P}(0,\lambda,\xi+\lambda
t)\exp\left({-\frac12\sigma^2_1|\lambda|^2t-\frac{\sigma^2_2
t}{6}\left(|\lambda|^2 t^2+3 t
(\lambda,\xi)+3|\xi|^2\right)}\right).
\end{equation}
The inverse Fourier transform (in the distributional sense) allows
to find the density $P(t,x,u)$:
$$P(t,x,u)=\dfrac1{(2\pi)^{2n}}\int\limits_{\mathbb{R}^n}
\int\limits_{\mathbb{R}^n}e^{i(\lambda,x)}e^{i(\xi,u)}\tilde{P}(t,\lambda,\xi)\,d\lambda
d\xi=$$
$$=\dfrac1{(2\pi)^{2n}}\int\limits_{\mathbb{R}^n}\rho_0(s)\int\limits_{\mathbb{R}^n}
\int\limits_{\mathbb{R}^n}e^{-\frac12\Lambda^2}
e^{-\frac12\Xi^2}e^{-\frac{|u_0(s)-u|^2}{2t\sigma_1^2}-\frac{6|(u_0(s)+u)\frac
t2+s-x|^2}{12\sigma_1^2t+\sigma_2^2t^3}}d\lambda d\xi ds=$$
\begin{equation}
\label{s_plotn}=C\int\limits_{\mathbb{R}^n}\rho_0(s)\,e^{-\frac{|u_0(s)-u|^2}{2t\sigma_1^2}-\frac{6|(u_0(s)+u)\frac
t2+s-x|^2}{12\sigma_1^2t+\sigma_2^2t^3}}ds,
\end{equation}
where
$$\Lambda=\Lambda(t,\lambda,s)=\dfrac{\sqrt{12\sigma_1^2t+\sigma_2^2t^3}}{2\sqrt{3}}\lambda+\dfrac{2i\sqrt{3}}{\sqrt{12\sigma_1^2t+\sigma_2^2t^3}}((u_0(s)+u)\frac
t2+s-x),$$
$$\Xi= \Xi(t,\xi,s)=\sigma_2\sqrt{t}\xi+\dfrac1{\sigma_2\sqrt{t}}\left(\dfrac{\sigma_2^2}2\lambda
t^2+iu_0(s)-u\right), \quad C=\left(\frac{\sqrt{3}}{\pi\sigma_2
t\,\sqrt{ 12\sigma_1^2+\sigma_2^2t^2}}\right)^{n}.$$

From (\ref{s_plotn}) we have
\begin{equation} \label{plotn}
\rho_\sigma(t,x)=\left(\dfrac{\sqrt{3}}{(\sqrt{2\pi
t(3\sigma_1^2+\sigma_2^2t^2)}}\right)^n\int\limits_{\mathbb{R}^n}\rho_0(s)e^{-\frac{3|u_0(s)t+s-x|^2}{2t(3\sigma_1^2+\sigma_2^2t^2)}}ds,
\end{equation}
\begin{equation}
\label{sol_u_sdu}
u_\sigma(t,x)=\dfrac1{2(3\sigma_1^2+\sigma_2^2t^2)}\dfrac{\int\limits_{\mathbb{R}^n}F(t,x,s)\rho_0(s)
e^{-\frac{3|u_0(s)t+s-x|^2}{2t(3\sigma_1^2+\sigma_2^2t^2)}}ds}{\int\limits_{\mathbb{R}^n}\rho_0(s)
e^{-\frac{3|u_0(s)t+s-x|^2}{2t(3\sigma_1^2+\sigma_2^2t^2)}}ds},
\end{equation}
where $F(t,x,s)=6\sigma_1^2u_0(s)-\sigma_2^2t(tu_0(s)+3(s-x))$.

\section{Properties of the ``perturbed velocity''}

\begin{prop}\label{prop1}
Let $u_0$ and $\rho_0>0 $ be  functions of class
$C^1(\mathbb{R}^n)\cap C_b(\mathbb{R}^n).$  If $ t_*(u_0)>0$ is a
moment of time such that the solution to the Cauchy problem
(\ref{equ_Burg}) with the initial condition $u_0$ keeps this
smoothness for $0<t<t_*(u_0)\le+\infty,$ then
$\hat{u}_\sigma(t,x)$ tends to a solution of problem
(\ref{equ_Burg}) as $|\sigma|\rightarrow 0$ for any fixed
$(t,x)\in\mathbb{R}^{n+1},\,0<t<t_*(u_0)$.
\end{prop}
\proof Let us denote by $J(u_0(x))$ the Jacobian matrix of the map
$\,x\longmapsto u_0(x).$ As it was shown in \cite{protter}\,(Theorem
1), if $J(u_0(x))$ has at least one eigenvalue which is negative for
a certain point $x\in{\mathbb R}^n,$  then the classical solution to
(\ref{equ_Burg}) fails to exist beyond a positive time $t_*(u_0).$
Otherwise, $t_*(u_0)=\infty.$ The matrix $C(t,x) = (I+t J(u_0(x))),$
where $\,I\,$ is the identity matrix, fails to be invertible for
$t=t_*(u_0).$

From (\ref{sol_u_sdu}) and (\ref{plotn}) we have
$$u_\sigma(t,x)\rho_\sigma(t,x)\,=\,(u^I_\sigma(t,x)-u^{II}_\sigma(t,x))\rho_\sigma(t,x),$$
where
$$u^I_\sigma(t,x)\,=\,
K \,\int\limits_{\mathbb{R}^n}u_0(s)\rho_0(s)
e^{-\frac{3|u_0(s)t+s-x|^2}{2t(3\sigma_1^2+\sigma_2^2t^2)}}ds,$$
$$u^{II}_\sigma(t,x)=\frac{3K\sigma_2^2t}{2(3\sigma_1^2+\sigma_2^2t^2)}\,\int\limits_{\mathbb{R}^n}(tu_0(s)+s-x)
\rho_0(s)e^{-\frac{3|u_0(s)t+s-x|^2}{2t(3\sigma_1^2+\sigma_2^2t^2)}}ds,$$
$K=\dfrac{\sqrt{3}}{\sqrt{2\pi t(3\sigma_1^2+\sigma_2^2t^2)}}.$
Using the weak convergence of measures and the fact that $\rho_0$
and $u_0$ are continuous and bounded, we get for $u^I_\sigma(t,x)$:
$$\lim\limits_{|\sigma|\rightarrow
0}u_\sigma^I(t,x)\,=\,\dfrac{\int\limits_{\mathbb{R}^n}u_0(s)\rho_0(s)\lim\limits_{\substack{\sigma_1\rightarrow
0\\
\sigma_2\rightarrow 0}}\left(\frac{\sqrt{3}}{\sqrt{2\pi
t(3\sigma_1^2+\sigma_2^2t^2)}}\right)^ne^{-\frac{3|u_0(s)t+s-x|^2}{2t(3\sigma_1^2+\sigma_2^2t^2)}}ds}{\int\limits_{\mathbb{R}^n}\rho_0(s)\lim\limits_{\substack{\sigma_1\rightarrow
0\\
\sigma_2\rightarrow 0}}\left(\frac{\sqrt{3}}{\sqrt{2\pi
t(3\sigma_1^2+\sigma_2^2t^2)}}\right)^ne^{-\frac{3|u_0(s)t+s-x|^2}{2t(3\sigma_1^2+\sigma_2^2t^2)}}ds}=$$
$$
=\,\dfrac{\int\limits_{\mathbb{R}^n}u_0(s)\rho_0(s)\delta_{p(t,x,s)}ds}{\int\limits_{\mathbb{R}^n}\rho_0(s)
\delta_{p(t,x,s)}ds},
$$
with $p(t,x,s) = u_0(s)t+s-x\,,$ where $\delta_y$ is the Dirac
measure as $y\in {\mathbb R}^n.$ We can then on the basis on the
invertibility of $C(t,x)$ use locally the implicit function
theorem and find $s=s_{t,x}(p).$

Therefore,
$$\lim\limits_{|\sigma|\rightarrow
0}u^I_\sigma(t,x)\,=\,
\dfrac{\int\limits_{\mathbb{R}^n}u_0(s_{t,x}(p))\rho_0(s_{t,x}(p))\,\det
(C(t,s_{t,x}(p)))^{-1}\,\delta_p\,(ds_{t,x})}
{\int\limits_{\mathbb{R}^n}\rho_0(s_{t,x}(p))\,\det(C(t,s_{t,x}(p)))^{-1}\,\delta_p\,(ds_{t,x})\,}=$$
$$=\,u_0(s_{t,x}(0)).$$
Now we prove that $\lim\limits_{|\sigma|\rightarrow
0}u^{II}_\sigma(t,x)\,=\,0$. We change variable $\tilde
s=\frac{\sqrt{3}(u_0(s)t+s-x)}{\sqrt{t(3\sigma_1^2+\sigma_2^2t^2)}}$
and obtain
$$\lim\limits_{|\sigma|\rightarrow
0}u^{II}_\sigma(t,x)\,=\,\lim\limits_{|\sigma|\rightarrow
0}\dfrac{\int\limits_{\mathbb{R}^n}\sigma_2^2t\rho_0(\tilde
s)\exp{(-\tilde s^2/2)d\tilde s}}{\rho_\sigma(t,x)}=0.$$

Let us introduce the new notation  $s_0(t,x)\equiv s_{t,x}(0).$
Then the following vectorial equation holds:
\begin{equation}
\label{usl}u_0(s_0(t,x))t+s_0(t,x)-x=0.
\end{equation}
The function $u(t,x)=u_0(s_0(t,x))$ satisfies the Burgers equation.
To check this fact it is enough  to differentiate (\ref{usl}) with
respect to $t$ and $x_j$ and find $s_{0,t}$ and $s_{0,x}$. The more
detailed proof you can find in \cite{AKR} in Proposition 1.
$\,\square$

It is important to note that $s_0(t,x)$ is unique for all $t$ for
which the solution to the Burgers equation $u(t,x)$ is smooth.

\begin{rem}Proposition \ref{prop1} can naturally be  extended to
the class of functions $\rho_0$ such that there exists a sequence
$\rho_0^\varepsilon\in C^1(\mathbb{R}^n)\cap C_b(\mathbb{R}^n)$
converging to $\rho_0$ as $\varepsilon\rightarrow 0$ almost
everywhere. In this case $u_\sigma(t,x)$ tends to a solution of
problem (\ref{equ_Burg}) as $\sigma\rightarrow 0$ almost
everywhere on $(t,x)\in\mathbb{R}^{n+1},\,0<t<t_*(u_0)$.
\end{rem}

\begin{thm}
The scalar function $\rho_{\sigma}(t,x)$ and the vector-function
$u_{\sigma}(t,x),$ defined in (\ref{u_sdu}), solve the following
system:
\begin{equation}
\label{sist_obw1}\dfrac{\partial\rho_\sigma}{\partial t}\,+\,div_x
(\rho_{\sigma}u_{\sigma})\,=\,\dfrac12\sigma_1^2\sum\limits_{k=1}^{n}\dfrac{\partial^2\rho_\sigma}{\partial
x_k^2},
\end{equation}
\begin{equation}
\label{sist_obw2}\dfrac{\partial(\rho_{\sigma}u_{\sigma,i})}{\partial
t}\,+\,
div_x(\rho_{\sigma}u_{\sigma,i}\,u_{\sigma})\,=\,\dfrac12\sigma_1^2\sum\limits_{k=1}^{n}\dfrac{\partial^2(\rho_{\sigma}u_{\sigma,i}))}{\partial
x_k^2}\,-I_{\sigma}
\end{equation}
where $i=1,..,n,\,t\ge 0$ and
$I_{\sigma}=\int\limits_{\mathbb{R}^n}(u_{i}-u_{\sigma,i})((u-u_{\sigma}),\nabla_x
P(t,x,u))du$.
\end{thm}

\proof Integrating  (\ref{Fok-Plank}) with respect of $u$ we get:
$$\dfrac{\partial\rho_\sigma}{\partial t}\,+\,div_x
(\rho_{\sigma}u_{\sigma})\,=\,\dfrac12\sigma_1^2\sum\limits_{k=1}^{n}\dfrac{\partial^2\rho_\sigma}{\partial
x_k^2}+\dfrac12\sigma_2^2\sum\limits_{k=1}^{n}\int\limits_{\mathbb{R}^n}\dfrac{\partial^2P(t,x,u)}{\partial
u_k^2}\,du.$$ The last integral is equal to zero, because $P(t,x,u)$
and all its derivatives tend to zero as $|u|\to\infty$.  Thus, we
obtain (\ref{sist_obw1}).

To prove (\ref{sist_obw2}) we note that the definitions of
$u_\sigma(t,x)$ and $\rho_\sigma(t,x)$ imply
\begin{center}
$$\dfrac{\partial(\rho_{\sigma}u_{\sigma,i})}{\partial
t} =\dfrac{\partial}{\partial
t}\int\limits_{\mathbb{R}^n}u_iP\,du=\int\limits_{\mathbb{R}^n}u_iP_t\,du
=-\int\limits_{\mathbb{R}^n}u_i(u,\nabla_x P)du+$$
\begin{equation}
\label{p2_1}+\dfrac12\sigma^2_1\sum\limits_{k=1}^{n}\int\limits_{\mathbb{R}^n}u_i\dfrac{\partial^2
P}{\partial
x_k^2}\,du+\dfrac12\sigma^2_2\sum\limits_{k=1}^{n}\int\limits_{\mathbb{R}^n}u_i\dfrac{\partial^2P}{\partial
u_k^2}\,du,
\end{equation}
\end{center}
where $P_t\equiv \frac{\partial }{\partial t} P.$ Early we noted
that the last integral is equal to zero.

Further, for $i=1,..,n$ we have
$$div_x(\rho\,u_{\sigma,i}\,u_\sigma)=
u_{\sigma,i}div_x(\rho_\sigma
u_\sigma)+u_\sigma\partial_x(\rho_\sigma u_{\sigma,i})-u_\sigma
u_{\sigma,i}\partial_x\rho_\sigma=$$
\begin{equation}
\label{p2_2}=\int\limits_{\mathbb{R}^n}u_{\sigma,i}(u,\nabla_xP)\,du+\int\limits_{\mathbb{R}^n}u_i(u_\sigma,\nabla_xP)\,du-\int\limits_{\mathbb{R}^n}u_{\sigma,i}(u_\sigma,\nabla_xP)\,du,
\end{equation}
where $i,k=1,...,n$. Equation (\ref{sist_obw2}) follows immediately
from (\ref{p2_1}) and (\ref{p2_2}). $\,\square$

Let us set $\rho(t,x)=\lim\limits_{|\sigma|\rightarrow
0}\rho_\sigma(t,x)$ and $\bar
u(t,x)=\lim\limits_{|\sigma|\rightarrow 0} u_\sigma(t,x).$
\begin{thm}\label{prop3}
Assume that $(\rho(t,x),\bar u(t,x)),$ the limits  of
$(\rho_\sigma,u_\sigma)$ as $|\sigma|\to 0,$ are $C^1$ -- smooth
bounded functions for $(t,x)\in \Omega:=[0,t_*(u_0))\times {\mathbb
R}^n,\, t_*(u_0)\le \infty.$ Then they solve in $\Omega$ the
pressureless gas dynamics system
\begin{equation}\label{sist_pred12}\partial_t
\rho+div_x(\rho \bar u)=0,\quad
\partial_t (\rho \bar u)+\nabla_x(\rho \bar u \otimes \bar u)=0.\end{equation}
\end{thm}

\proof As follows from Proposition \ref{prop1}, the function $\bar
u(t,x)$ is a $C^1$ - solution of the non-viscous Burgers equation.
Further,  (\ref{sist_obw1}) is a linear parabolic equation with
respect to $\rho_\sigma,$ hence the limit as $|\sigma|\to 0$ reduces
it to the first equation in (\ref{sist_pred12}) (continuity
equation). The second equation in (\ref{sist_pred12}) is a corollary
of the non-viscous Burgers equation and continuity equation for
smooth solutions. $\square$

\begin{rem} Proposition \ref{prop3} implies that the integral
term on the right-hand side of (\ref{sist_obw2}) vanishes as
$|\sigma|\to 0$ in the case of smooth limit functions $\rho$ and
$\bar u.$
\end{rem}

\section{A regularization of the pressureless gas dynamics}
After the moment of shock formation in the solution to the Burgers
equation the limit  of \eqref{sist_obw1}-\eqref{sist_obw2} as
$|\sigma|\to 0$ is not equivalent to the pressureless gas dynamics
system \eqref{sist_pred12}, the integral term in \eqref{sist_obw2}
does not vanish and can be interpreted as a gradient of special
pressure (spurious pressure). It is well known that for the
pressureless gas dynamic system after the moment of singularity
formation there arise a question how to define a generalized
solution, and the answer depends on the prescribed type of
interaction between particles. For the sticky particles model, where
two particles move together after interaction, one need to introduce
a strongly singular solution with a delta-singularity in the density
component. If we assume that particles do not interact, we get the
generalized solution in the sense of free particles, whose
singularities are only shocks. In this sense we can call the  limit
of system \eqref{sist_obw1}-\eqref{sist_obw2} as $|\sigma|\to 0$  a
regularization of the pressureless gas dynamics equations. There
exist an algorithm how to pass from the model of free particles to
the model of sticky particles (see \cite{AKR} for details in the
case $\sigma_2=0$).

For discontinuous initial data analogously to \cite{AKR} we
introduce a generalized solution in the sense of free particles:
\begin{defn}
We call the couple of functions $(\rho_{FP}(t,x),u_{FP}(t,x))$ the
generalized solution to the Cauchy problem for \eqref{sist_pred12}
in the sense of \,free \,particles \,(FP-generalized \,solution) \,
subject \,to \,initial \,data \\$(\rho_0(x),u_0(x))$ $ \in {\mathbb
L}^2_{loc}({\mathbb R}^n)\cap {\mathbb L}_{\infty}({\mathbb R}^n)$,
if for almost all $(t,x)\in \mathbb{R}_+\times \mathbb{R}^n$
$$\rho_{FP}(t,x)=\lim\limits_{\varepsilon\rightarrow
0}(\lim\limits_{\sigma\rightarrow
0}\rho_\sigma^{\varepsilon}(t,x)),\qquad
u_{FP}(t,x)=\lim\limits_{\varepsilon\rightarrow
0}(\lim\limits_{\sigma\rightarrow
0}{u}_\sigma^{\varepsilon}(t,x)),$$ where
$(\rho_\sigma^{\varepsilon}(t,x), {u}_\sigma^{\varepsilon}(t,x))$
correspond to  initial data $\, \rho_0^{\varepsilon}\,=\,
\eta_\varepsilon *\rho_0,\quad
u_0^{\varepsilon}\,=\,\eta_\varepsilon *u_0\,\,,$ where
$\eta_\varepsilon(x)$ is the standard averaging kernel.
\end{defn}
In \cite{AKR} for $\sigma_2=0$ we solve the Riemann problem in the
sense of free particles, this can be done analogously in the general
case $\sigma_2>0$.

It is very important that the double limit procedure used in the
definition of the FP-solution plays role only for the case of a
central rarefaction wave, where it helps to find a unique stable
solution in the vacuum domain (we refer to \cite{AKR} again). For
the compression wave the FP-solution can be computed directly by
formulae \eqref{plotn}, \eqref{sol_u_sdu} using only one limit pass
$|\sigma|\to 0$. Thus, after the moment of the shock formation
$t_*(u_0)$ in the Burgers equation the solution to the respective
conservation law
\begin{equation}\label{sist_press}\partial_t
\rho+div_x(\rho \bar u)=0,\quad
\partial_t (\rho \bar u)+\nabla_x(\rho \bar u \otimes \bar u)+ \nabla p=0,\quad \nabla p =
\lim\limits_{|\sigma|\to 0}I_\sigma\end{equation} can be found from
\eqref{plotn}, \eqref{sol_u_sdu} as $|\sigma|\to 0$. By definition
the solution to \eqref{sist_press} is the FP-solution to
\eqref{sist_pred12}, thus \eqref{sist_press} can be considered as a
regularization of \eqref{sist_pred12}.

\section{The 1D model for polytropic gas}

Let us {\it assume} that the gas that obeys system
\eqref{sist_press} is polytropic, i.e. the state equation is
$p=A\rho^\gamma$, $\gamma$ is the adiabatic exponent, $A$ is a
positive constant. Therefore we have to supplement
\eqref{sist_press} with equation $$\partial_tp+u\partial_xp+\gamma
p\partial_xu=0.$$ We consider the 1D Riemann problem in the
compression case $u_2<0$: $\rho_0(x)=\rho_1+\rho_2\theta(x)$,
$u_0(x)=u_1+u_2\theta(x)$, where $\rho_i$, $u_i$, $i=1,2$ are
constants, In \cite{AKR} the solution was obtained (for
$\sigma_2=0$, the general case is analogous):
$$\rho(t,x)=\rho_1+(\rho_1+\rho_2)\theta(x-u_1t)-\rho_1\theta(x-(u_1+u_2)t),$$
$$u(t,x)=u_1+\dfrac{\rho_1+\rho_2}{2\rho_1+\rho_2}u_2\theta(x-u_1t)-\dfrac{\rho_1}{2\rho_1+\rho_2}u_2\theta(x-(u_1+u_2)t),$$
$$p(t,x)=\dfrac{\rho_1(\rho_1+\rho_2)}{2\rho_1+\rho_2}u_2^2(\theta(x-(u_1+u_2)t)-\theta(x-u_1t)).$$
\begin{prop} In assumption that the state equation is polytropic, the adiabatic exponent $\gamma\to 3$ as  the ratio
$\dfrac{\rho_1}{\rho_2}\to 0$.
\end{prop}
\proof We have three Hugoniot conditions, that follow from the mass
conservation, momentum conservation (\eqref{sist_press}) and the
conservation of total energy $E=\dfrac12\rho
u^2+\dfrac{p}{\gamma-1}$. We denote by $f_-$ and $f_+$ the values of
$f(x)$ before and after shock, respectively. For the shock moving
with the speed $D=u_1+u_2$ we have $\rho_+=2\rho_1+\rho_2,$
$\rho_-=\rho_1$, $u_+=u_1+\dfrac{\rho_1+\rho_2}{2\rho_1+\rho_2}u_2$,
$u_-=u_1$, $p_+=\dfrac{\rho_1(\rho_1+\rho_2)}{2\rho_1+\rho_2}u_2^2$,
$p_-=0$. The third  Hugoniot condition gives
$D=\dfrac{[Eu+pu]}{[E]}$ and
$\gamma=\dfrac{3\rho_1+\rho_2}{\rho_1+\rho_2}.$ This implies the
proposition. The consideration of the second shock gives the same
result. $\square$

The adiabatic exponent $\gamma=1+\frac{2}{n}$, where $n$ is the
space dimension, corresponds to the monoatomic gas.

\end{document}